\documentclass[a4paper]{article}
\usepackage{hyperref,color,ifthen}
\usepackage{amsmath,amsfonts,amssymb}
\usepackage{graphicx}

\newtheorem{Example}{Example}[section]

\newtheorem{Statement}{Statement}[section]
\newenvironment{Proof}{{\it Proof}\ }{$\Box$}

\def\tx#1{{\rm #1}\,}

\title{Gauss-Bonnet type theorems in any dimension
\footnote{This work was partially supported by RFBR 01-01-00546 grant.}}
\author{Valery DOLOTIN\footnote{e-mail:\ {\it vd@gate.itep.ru}}\\
Institute of Theoretical \& Experimental Physics, Moscow}

\begin{document}
\maketitle
\vspace{-9cm}
\begin{flushright}
Preprint ITEP-TH-88/01
\end{flushright}
\vspace{7.5cm}
\begin{abstract}
Given a construction of smooth homotopy class invariants of smooth
immersions $M^n\to{\bf R}^{n+k}$. The particular case of $k=1,\
n\ge 1$ is a sequence of non-zero integrals, where the $n=2$ term
is the Gauss-Bonnet integral
\end{abstract}

\section{Basic example}

Take a smooth curve $\gamma=(x(t),y(t))\subset{\bf R}^2$.

\begin{Statement}
$$I(\gamma):=\int\frac{\dot y\ddot x-\dot x\ddot y}{\dot x^2+\dot y^2}dt=2\pi$$
\end{Statement}
\begin{Proof}
Denote $u:=\dot x,\ v:=\dot y$. Then the integral becomes
$$\int_\gamma\frac{vdu-udv}{u^2+v^2}$$
where $\frac{vdu-udv}{u^2+v^2}$ is a closed from on ${\bf R}^2-0$
being a generator of $H^1({\bf R}^2-0;{\bf R})$.
\end{Proof}

Note, that the above statement is the smooth analogue of the
Euclidean geometry theorem about the sum of exterior angles of a
polygon.

If the curve has self-intersections, then in general
$I(\gamma)=2\pi k$ where $k$ is the "winding number" of the
tangent vector while moving along the curve. If we take a homotopy
connecting two curves with different winding numbers then we have
to pass trough a curve with a singularity of the form
$\frac{vdu-udv}{u^2+v^2}$ at some point. Those singular curves
make a discriminant set in the space of smooth curves.

\section{General construction}

\subsection{Codimension $1$ case}

For each space $P=(p_0,\dots,p_n)$ we have a canonical closed form
$$\omega:=\frac{\sum_i (-1)^ip_idp_0\wedge\dots\widehat{dp_i}\dots\wedge dp_n}{(p_0^2+\dots+p_n^2)^{(n+1)/2}}$$
on $P-0$.

Then for a smooth immersion of $n$-manifold $M\to{\bf R}^{n+1}$ in
each local chart $(t_1,\dots,t_n)$ on $M$ we have a Jacobian
$n\times(n+1)$ matrix
$$J:=\left(\frac{\partial x_i}{\partial t_j}\right)$$
If we substitute via Plucker embedding
$p_i=\Delta_{i}:=\Delta_{0\dots\widehat i\dots n}(J)$ (being the
minor build on all except $i$-th columns of $J$) into $\omega$
then we get an integration form $\omega(\Delta)$ on $M$. This form
is correctly defined since the orientation preserving change of
coordinates between charts multiplies all $\Delta_{\widehat i}$ by
the same number (the Jacobian of the change) and $\omega$ is
invariant under such multiplication. For a non-orientable manifold
the form $\omega(\Delta)$ may be interpreted as a two-valued from,
or a form on a double covering of $M$.

Since $\omega$ is closed then

\begin{Statement}
$$\int_M\omega(\Delta)=2^n\pi k$$
where $k$ is an integer invariant of a smooth homotopy class of
smooth immersions of $M$.
\end{Statement}

Note that in codimension 1 case this invariant equals to the
degree of the corresponding Gauss map $M^n\to S^n$.

\subsubsection{Example: Gauss-Bonnet theorem}

Here $n=2$ and
$$\omega=\frac{p_0dp_1\wedge dp_2-p_1dp_0\wedge dp_2+p_2dp_0\wedge dp_1}{(p_0^2+p_1^2+p_2^2)^{3/2}}$$

For a smooth immersion $M^2\to{\bf R}^3,\ (t_1,t_2)\mapsto
(x_1,x_2,x_3)$ with the Jacobian
$$J=\left(\frac{\partial x_i}{\partial t_j}\right)=:\left(x_{i,j}\right)$$
and $p_0=\Delta_0(J)=x_{1,1}x_{2,2}-x_{2,1}x_{1,2},\
p_1=\Delta_1(J)=x_{0,1}x_{2,2}-x_{2,1}x_{0,2},\
p_2=\Delta_2(J)=x_{0,1}x_{1,2}-x_{1,1}x_{0,2}$ we have the
integral
$$\int_M\omega(\Delta)=\int_M\frac{\sum_{i,j}\Delta_i\Delta_jdx_{i,1}\wedge dx_{j,2}}
{(\sum\Delta_i^2)^{3/2}}=
\int\frac{\sum_{i,j}\Delta_i\Delta_j(x_{i,11}x_{j,22}-x_{i,12}x_{j,21})dt_1dt_2}{(\sum\Delta_i^2)^{3/2}}$$

One can check that $\sum\Delta_i^2=EG-F^2$ is the determinant (in
Gauss notations) of the metric form induced on $M$ from the
Euclidean metric on ${\bf R}^3$. Now for the numerator
$$\sum_{i,j}\Delta_i\Delta_j(x_{i,11}x_{j,22}-x_{i,12}x_{j,21})=$$
$$=\left|
\begin{array}{ccc}
x_{1,11}&x_{2,11}&x_{3,11}\\
x_{1,1}&x_{2,1}&x_{3,1}\\
x_{1,2}&x_{2,2}&x_{3,2}
\end{array}
\right| \left|
\begin{array}{ccc}
x_{1,22}&x_{2,22}&x_{3,22}\\
x_{1,1}&x_{2,1}&x_{3,1}\\
x_{1,2}&x_{2,2}&x_{3,2}
\end{array}
\right|$$$$- \left|
\begin{array}{ccc}
x_{1,12}&x_{2,12}&x_{3,12}\\
x_{1,1}&x_{2,1}&x_{3,1}\\
x_{1,2}&x_{2,2}&x_{3,2}
\end{array}
\right| \left|
\begin{array}{ccc}
x_{1,21}&x_{2,21}&x_{3,21}\\
x_{1,1}&x_{2,1}&x_{3,1}\\
x_{1,2}&x_{2,2}&x_{3,2}
\end{array}
\right|=ln-m^2$$ is the determinant of the second fundamental form
(in Gauss notations). Since the area element on $M$ is
$dS=\sqrt{EG-F^2}dt_1dt_2$ then we have
$$\int_M\omega(\Delta)=\int\frac{ln-m^2}{(EG-F^2)^2}dS=2\pi\chi(M)$$
which is the Gauss-Bonnet formula.

\subsection{Affine Gauss map}

For a smooth immersion $\phi:M^n\to {\bf R}^{n+k}$ taking the
minors of its Jacobian we get a map ${\cal J}:M^n\to G_{n,n+k}$,
%where $\tilde G_{n,n+k}$ is the space of oriented $n$-subspaces in
%${\bf R}^{n+k}$ (which is a double covering of the Grassmannian
%$G_{n,n+k}$). This
which is the higher codimension analogue of the Gauss map.

Then for any closed $n$-form $\omega\in H^n(G_{n,n+k})$ its
pullback ${\cal J}^*\omega$ gives us a smooth homotopy class
invariant $\int_M{\cal J}^*\omega$ of smooth immersions of $M$.

The generators of $H^n(G_{n,n+k})$ have a natural representation
in terms of Plucker coordinates $p_I$ on $G_{n,n+k}$ as
$$\omega=\frac{\sum \varphi_{I_1\dots I_n}(p)dp_{I_1}\wedge\dots\wedge dp_{I_n}}{|p|^d}$$
where $|p|:=(\sum p_I^2)^{1/2}$, and $\varphi_{I_1\dots I_n}(p)$
are homogeneous polynomials of degree $d-n$.

\subsubsection{Example: $n=2$, Kahler form}

On the space ${\bf R}^{2N}-\{0\}$ we have a closed 2-form
$$\omega=
%\sum_{1\le i<j\le d}
(\sum_{i,j}(p_idp_i+p_{i+N}dp_{i+N})\wedge(p_{j+N}dp_j-p_jdp_{j+N})$$
$$-(p_jdp_i+p_{j+N}dp_{i+N})\wedge(p_{j+N}dp_i-p_jdp_{i+N}))\,/\,
|p|^4$$
%{(p_idp_i+p_{i+d}dp_{i+d})\wedge((p_j+p_{j+d})dp_j-(p_j-p_{j+d})dp_{j+d})}
%{(\sum_1^{2d}p_i^2)^2}$$
%p_i(p_j+p_{j+d})dp_i\wedge dp_{j}+p_{i+d}(p_j+p_{j+d})dp_{i+d}\wedge
%dp_{j}+$$ $$p_i(p_{j+d}-p_{j})dp_{i}\wedge
%dp_{j+d}+p_{i+d}(p_{j+d}-p_{j})dp_{i+d}\wedge dp_{j+d}
%$$\omega_i=\sum_{1\le i<j\le d}\frac
%{(p_{i+d}dp_i-p_idp_{i+d})\wedge((p_j+p_{j+d})dp_j+(p_j-p_{j+d})dp_{j+d})}
%{(\sum_1^{2d}p_i^2)^2}$$
%$^{\left(^2_4\right)}$
which is the imaginary (non-zero) part of the Kahler form
$\partial\bar\partial\log(|z_1|^2+\dots+|z_N|^2)$ written in terms
of complex coordinates $z_j=p_j+ip_{j+N}$.

Note, that as a form on ${\bf R}^{2N}-\{0\}$ it is exact, but it
is invariant under the natural ${\bf C}^*$-action on ${\bf
R}^{2N}-\{0\}$ and its pull-back to the space ${\bf C}P^{N-1}$ of
${\bf C}^*$-orbits gives a non-trivial cohomology class.

In particular, for a smooth complex submanifold $X^n\subset {\bf
C}^{n+k}$ the corresponding Jacobian matrix, combined with the
Plucker embedding gives us a map $X^n\to{\bf C}P^{\left(^{\
n}_{n+k}\right)-1}$. Then the pull-back of $\omega$ is known as
the first Chern class $c_1(X\cap H)$ of hyperplane section line
bundle over $X$. On the other hand, taking $X^n$ as a real
$2n$-dimensional immersion into ${\bf R}^{2(n+k)}$ we have the
pull-back of $\omega$ in terms of real partial derivatives
$\partial x_i/\partial t_j,\ i=1,\dots,2(n+k),\ j=1,\dots,2n$.
Since the real $d$-closeness of $\omega$ is just the algebraic
property of its expression in terms of these real partial
derivatives, then $\int_M J^*\omega^n$ is a differential invariant
for any immersion $M^{2n}\to{\bf R}^{2(n+k)}$.
%
%Then for an immersion $M^2\to{\bf R}^{2+k}$, with the
%corresponding "Gauss" map $\gamma:M^2\to \tilde G_{2,2+k}$ and
%Kahler forms $(\omega_r,\omega_i)$ on the Plucker space ${\bf
%R}^{\left(^{\ 2}_{2+k}\right)}-0\hookleftarrow\tilde G_{2,2+k}$,
%we get two differential invariants on $M$:
%$$\int_M\gamma^*\omega_r\quad{\rm and}\quad \int_M\gamma^*\omega_i\ .$$

\subsection{Projective Gauss map}

For each point $x$ of a smooth projective manifold $M^n\subset{\bf
R}P^N$ we have a $(n+1)$-dimensional subspace in ${\bf R}^{N+1}$,
tangent to $M$ at $x$, given by the Jacobian matrix of the cone
over $M$ in ${\bf R}^{N+1}$. So we have a (projective) Gauss map
$J:M^n\to G_{n+1,N+1}$.

\begin{Example}
Take $n=1,\ N=2$, i.e. a curve $\gamma\subset {\bf R}P^2$. Take a
chart $U_0:=\{(x_0,x_1,x_2)|\ x_0\ne 0\}\subset{\bf R}^3-\{0\}$
and write the 2-parametric equation of the cone over $\gamma$ as
$x_0=s,\ x_1=sf_1(t),\ x_2=sf_2(t)$. The corresponding Jacobian
matrix is
$$J=\left(
\begin{array}{ccc}

0 & sf_1' & sf_2'\\
1 & f_1  & f_2\\
\end{array} \right)$$

Take the form $\omega_0=\frac{p_2dp_1-p_1dp_2}{p_1^2+p_2^2}$ on
${\bf R}^3-\{0\}$. Then its pull-back $J^*\omega_0$ to the image
of Plucker embedding $p_i=\Delta_i(J)$ of the Gauss image of
$\gamma$ is the form $\frac{\dot f_2\ddot f_1-\dot f_1\ddot
f_2}{\dot f_1^2+\dot f_2^2}dt$, known from Section 1.

Note, that if $\gamma$ has a singularity in the $U_0$ chart then
$J^*\omega_0$ has a singularity at the corresponding point of
$\gamma$ (although the inverse is not true). Then we have a bundle
$\omega(\alpha):=\sum_0^2 \alpha_i\omega_i$ of singular forms on
${\bf R}^3$ of which the pull-backs to $\gamma$ give integrals,
being invariants if smooth homotopy classes of the Gauss map
$\gamma\to{\bf R}P^2-\{q_0,q_1,q_2\}$, where $q_i=\{p\in{\bf
R}P^2|\ \sum_{j\ne i}p_j^2=0 \}$.
\end{Example}

In the chart $U_i:=\{(x_0,\dots,x_N)|\ x_i\ne 0\}$ the cone over
$M$ may be given parametrically as $x_0=s\cdot
f_0(t_1,\dots,t_n),\dots,\ x_i=s,\dots,\ x_N=s\cdot
f_N(t_1,\dots,t_n)$ with the Jacobian
$$J=
\left(
\begin{array}{ccccc}

sf_{0,t_1} & \dots & 0      & \dots & sf_{N,t_1}\\
\vdots    &       & \vdots &       & \vdots\\
sf_{0,t_n} & \dots & 0      & \dots & sf_{N,t_n}\\
f_0       & \dots & \underbrace{1}_i      & \dots & f_N\\

\end{array}
\right)$$ where $f_{i,t_j}:=\frac{\partial f_i}{\partial t_j}$.
Having a singularity of $M$ in the chart $U_i$ implies setting
$\Delta_{ii_1\dots i_n}(J)=0$, for all $ i_1<\dots <i_n\subset
\{1,\dots,N \}$, which also implies the degeneracy of the whole
$J$. Then on the space of Plucker coordinates ${\bf
R}^{\left(^{\,n+1}_{N+1}\right)}$ any form
$$\omega=(\sum
\varphi_{I_1\dots I_n}(p)dp_{I_1}\wedge\dots\wedge
dp_{I_n})\,/\,|p|^d$$ (where $\varphi(p)$ are homogeneous of
degree $d-n$)
%$$\omega_i=(\sum
%\varphi_{I_1\dots I_n}(p)dp_{I_1}\wedge\dots\wedge
%dp_{I_n})\,/\,|p|_i^d$$ (where
%$$|p|_i=(\sum_{I=(ii_1\dots i_n)} p_I^2)^{1/2}$$ is the
%sum of squares of all Plucker coordinates containing fixed $i$ in
%the index, and $\varphi(p)$ are homogeneouse of degree $d-n$)
has a pull-back $J^*\omega$ to $M$ with singularities at all
singular points of $M$, regardless the chart of their location.

%Now take $\omega(\alpha):=\sum\alpha_i\omega_i$ with all
%$\omega_i$ being closed.

Now take $\omega$ to be closed.

\begin{Statement} For the pull-back $J^*\omega$ to $M$ the
integral
$$\int_MJ^*\omega$$ is a smooth homotopy class invariant of $M^n\subset
{\bf R}P^N$. \end{Statement}

Note, that even if $\omega$ is exact as a form on ${\bf
R}^{\left(^{\,n+1}_{N+1}\right)}-\{0\}$ its pull-back to the
projective manifold may give a non-trivial class (as it's the case
in the above example).

For a given immersion $M^n\to {\bf R}^{n+k}$ we may take a
composition of the affine Gauss map ${\cal J}_a:M^n\to
G_{k,n+k}\hookrightarrow {\bf R}P^{\left(^{\ k}_{n+k}\right)-1}$
with the projective one ${\cal J}_p:{\cal J}_a(M)\to
G_{n+1,\left(^{\ \, k}_{n+k}\right)}$ and use the corresponding
invariants.

\section{Gauss-Bonnet type invariants for polyhedra}

\noindent\underline{$n=1$}\ This is the above mentioned theorem on the sum of the exterior angles of a polygon.
%\linebreak

\noindent\underline{$n=2$\ }
For a vertex $\nu$ of a triangulation the sum of angles over all triangles $\tau$ of its star 
$$\sum_{\tau\in St_\nu }(\alpha_\tau+\beta_\tau+\gamma_\tau)=\pi|St_\nu|$$ 
(where $|St_\nu|$ is the number of triangles) is a constant, invariant with respect to different immersions of the star. If $\alpha_\tau$ denote angles adjesent to $\nu$ then we can rewrite this identity as 
$$2\pi-\sum_\tau\alpha_\tau=2\pi-\sum_\tau(\pi-(\beta_\tau+\gamma_{\tau+1}))$$
Note, that $\pi-(\beta_\tau+\gamma_{\tau+1})$ may by interpreted as a deformed exterior angle at a vertex $\mu\in Ln_\nu$ with ajasent $\beta_\tau$ and $\gamma_{\tau+1}$ which differs from the actual exterior angle as "much" as the immersion of the star at $\mu$ differs form lying in a $2$-plane. We call this difference by a curvature $K_1(\mu)$ of a $1$-dimensional link at the vertex $\mu$. Since $2\pi$ is the topological invariant of the link, which is equal to the sum of exterior angles when the link lies in a $2$-plane then we can rewrite the above equation as
$$E_2(\nu):=2\pi-\sum_{\tau\in St_\nu}\alpha_\tau=2\pi-\sum_{\mu\in Ln_\nu}K_1(\mu)$$ 
and call $E_2(\nu)$ the {\bf exterior $2$-angle} of $\nu$ via the immersion into $\mathbb{R}^3$. Now the analogue of the theorem on exterior angles of a polygon is 
\begin{Statement}
$$\sum_\nu E_2(\nu)=2^2\pi\kappa$$ is the topological invariant of the immersion.
\end{Statement}
%In particular for a sphere $\chi=1$.

\noindent\underline{Inductive step}

Having defined the exterior angle $E_k(\nu)$ for a vertex of a $k$-polyhedron in $\mathbb{R}^{k+1}$ and knowing the corresponding topological invariant $2^k\pi$ of a $k$-sphere we can define the curvature $K_k(\mu)$ at vertices of the $k$-link of a $k+1$-star immersed into $\mathbb{R}^{k+2}$. Then the $k+1$-exterior angle is defined by induction
$$E_{k+1}(\nu):=2^{k}\pi-\sum_{\mu\in Ln_\nu}K_k(\mu)$$
which may be regarded as a particular case of polyhedral Stokes formula, $2^{k}\pi-E_{k+1}(\nu)$ corresponding to integration form of degree $k+1$. 

The corresponding invariant of immersions into $\mathbb{R}^{k+2}$ is 
$$\sum_{\nu\in P}E_{k+1}(\nu)$$

%Take a vertex $\nu$ of a $2$-dimensional PL-manifold. For $\mu\in{\tx Ln}\nu$ take ${\tx %St }\nu\cap{\tx St}\mu$. Define an effective curvature of the link of $\nu$ at a given %vertex $\mu\in{\rm Ln}(\nu)$ as $2\pi-{\rm Ln}{\tx St}\nu\cap{\tx St}\mu$. To define the %curvature of a $2$-star we take 

\section{Discussion}

The above consideration implies, that having a differential
invariant of a manifold, written in terms of "intrinsic geometry"
of $M$ (such as Gauss-Bonnet), it may be convenient to represent
it as a restriction onto $M$ of some natural closed form defined
on the appropriate target space of a smooth immresion of $M$.

It would be interesting to see the polyhedral interpretation of
the integral $\int_M\omega(\Delta)$ for $\tx{codim} M\ge 2$ analogous to
that in the curvature case above. Provided having the polyhedral presentation of $\omega(\Delta)$ we could compute the value of the corresponding polyhedral form at vertices of any polyhedra, including non-smoothable ones. Since each polyhedron may be smoothed almost everywhere then each embedding of such almost everywhere smooth manifold has irremovable singularities at certain points corresponding to "non-smoothable" vertices of our triangulation. Then the value of the polyhedral form at non-smoothable vertices should distinguish them among others, giving some maximal possible value (like $2\pi$ value for the curvature form at the cusp singularity of a plane curve).

\end{document}